\documentclass{article}

 \usepackage{amsxtra,amsthm}
\usepackage{graphicx}
\usepackage{delarray}
\usepackage{mathptmx}
\usepackage{amsmath, amstext,amssymb,amsfonts}

\DeclareFontFamily{U}{mathx}{\hyphenchar\font45}
\DeclareFontShape{U}{mathx}{m}{n}{
      <5> <6> <7> <8> <9> <10>
      <10.95> <12> <14.4> <17.28> <20.74> <24.88>
      mathx10
      }{}
\DeclareSymbolFont{mathx}{U}{mathx}{m}{n}
\DeclareFontSubstitution{U}{mathx}{m}{n}
\DeclareMathAccent{\widecheck}{0}{mathx}{"71}
\DeclareMathAccent{\wideparen}{0}{mathx}{"75}

\begin{document}
{\centerline{\bf A note on $n$-divisible positive definite  functions}}

\vspace{6mm}
\centerline{{\large Saulius  Norvidas} }

\vspace{3mm}
{{\footnotesize Faculty of Mathematics and Informatics, Vilnius University, Akademijos str. 4, 08412 Vilnius, Lithuania;  \\
 ({\rm{e-mail: norvidas{@}gmail.com}})}}


\newtheorem{thm}{Theorem}
 \newtheorem{cor}[thm]{Corollary}
 \newtheorem{lem}[thm]{Lemma}
 \newtheorem{prop}[thm]{Proposition}
 \newtheorem{rem}[thm]{Remark}
  \newtheorem{exam}[thm]{Example}
  \newcommand{\R}{\mathbb R}    \newcommand{\Co}{\mathbb C}             \newcommand{\Z}{\mathbb Z}    \newcommand{\Zn}{\mathbb Z_n} \newcommand{\N}{\mathbb N} \newcommand{\M}{M({\mathbb R})}
  \newcommand{\PD}{PD(\mathbb R)}      \newcommand{\CH}{CH(\mathbb R)}      \newcommand{\Le}{L^1({\mathbb R})}  \newcommand{\B}{B(\R)} \newcommand{\A}{A(\R)} \newcommand{\PM}{M_{+}(\R)} \newcommand{\Pn}{PD_n(\mathbb R)} \newcommand{\Pe}{PD_{\infty}(\mathbb R)}\newcommand{\car}{\text{card}}
 \newcommand{\Dn}{D_n(f)}
\begin{abstract}
 Let $\PD$ be  the family of continuous  positive definite functions on $\R$. For an integer $n>1$, a $f\in\PD$ is called $n$-divisible if there is $g\in\PD$ such that $g^n=f$.  Some properties of infinite-divisible and $n$-divisible functions may differ in essence. Indeed, if $f$ is infinite-divisible, then for each integer $n>1$, there is an unique $g$ such that $g^n=f$, but there is a  $n$-divisible $f$ such that the factor $g$ in $g^n=f$ is generally not unique. In this paper, we discuss about how rich can be the class $\{g\in\PD: g^n=f\}$ for $n$-divisible $f\in\PD$ and obtain precise estimate for the cardinality of this class.

\end{abstract}


{{\it{ MSC}}} 2010: \  42A38 - 42A82 - 60E10
{\large{
\section{Introduction}
 We start with some notations and definitions. Let $\Z$, $\N$, $\R$ and $\Co$ be the families of  integers, positive integers,  real  and complex numbers, respectively. In the sequel, $\M$  denotes   the  Banach algebra of bounded regular complex-valued  Borel measures   on $\R$ with the convolution as multiplication. $\M$ is equipped with the usual total variation norm $\|\mu\|$ of  $\mu\in\M$.
 The Fourier-Stieltjes transform of  $\mu\in\M$ is given by
\[
\hat{\mu}(x)=\int_{\R}e^{-i x t} d\mu(t).
\]
A function $f : \R\to\Co$  is said to be positive definite if
\[
\sum_{j, k=1}^{m} f(x_j-x_k)c_j\overline{c}_k\ge 0
\]
 for each $m\in\N$ and all  $c_1,\dots, c_m\in \Co$,   $x_1, \dots, x_m \in\R$. Any such a function  satisfies
\begin{equation}
f(-x)=\overline{f(x)}
\end{equation}
for all  $x\in\R$.

Positive definite functions on groups  have a long history and have many applications in probability theory and such areas as stochastic processes [1]; harmonic analysis  [3],
potential theory [5] and spectral theory [2]. See [10] for other applications and details. The analysis of the properties of positive definite  functions
has  a vast literature, and the above list is only a small sample.

We will denote by $\PD$ the family of  continuous nontrivial   ($\not\equiv 0$) positive definite functions  on $\R$.   Note that if positive definite $f$ is continuous in a neighborhood of the origin, the it is uniformly continuous in $\R$ (see, e.g., [4, Corollary 1.4.10]). Bochner's  theorem gives a description of  $f\in\PD$  in terms of the Fourier transform. Namely, according to this theorem  (see, e.g., \cite[p. 71]{9}),  a continuous function $f:\R\to\Co$ is positive definite if and only if there exists  a nonnegative $\mu_f\in\M$ such that $f(x)=\widehat{\mu_f}(x)$, $x\in\R$. This statement implies, in particular, that
\begin{equation}
|f(x)|\le f(0)=\|\mu_f\|
\end{equation}
for all  $x\in\R$.
If,  in addition,   positive $\mu_f\in\M$  satisfies $\|\mu_f\|=1$, then on the language of probability theory, such  a $\mu_f$ and the function $f(x)=\widehat{\mu_f}(-x)$, $x\in\R$, are called a probability measure and its characteristic function, respectively  (see, e.g., \cite[p.p. 8-9]{17}).

Recall that in the probability theory a random variable $\xi$  is called $n$-divisible  for  certain $n\in\N$, $ n>1$,  if there exist independent and identically distributed random variables $\xi_1,\dots, \xi_n$ such that $\xi_1+\dots +\xi_n$ has the  sane distribution as $\xi$.    In terms of the characteristic function $f$ of a real-valued random variable, this means that  $f$ is $n$-divisible if there exists $g\in\PD$ such that
\begin{equation}
f(x)=  g^n(x)
\end{equation}
 for all  $x\in\R$.   Next, $f\in\PD$ is said to be infinite-divisible      if it is $n$-divisible for each $n\in\N$, $n\ge 2$.

 In the sequel, $\Pn$ and $\Pe$ denotes the families of $n$-divisible and infinite divisible functions in $\PD$, respectively. An early overview over divisibility of distributions is given in [12]. Until very recently, the vast majority of  divisible positive definite functions or divisible distributions  considered in the literature are also infinite-divisible. Important applications of $n$-divisibility is in modelling, for example of bug populations in entomology [7], or in financial aspects of various insurance models [14] and [15].

The motivation for our investigation comes: (i) partly from the fact that properties of functions in $\Pe$ has a rich literature (see, e.g.,   \cite{8} and \cite{13}), but the $n$-divisible functions have been studied much less; (ii)   partly from the fact that some properties of functions from  $\Pe$ and from $\Pn$  may differ in essence. One of those properties is the following: if $f\in\Pe$, then for each $n\in\N$, $n>1$, there is an unique $g\in\PD$ such that $g^n=f$, but there is  $n$-divisible $f$ such that the factor $g$ in $g^n=f$ is generally not unique. In this paper, we study the following problems:  (i) how rich can be the class $\{g\in\PD: g^n=f\}$; \ (ii) what properties of $f$  determine the size of  $\{g\in\PD: g^n=f\}$. We present several
precise estimates for the cardinality of this class. Also, the main results are validated via illustrative examples.

More precise, for $n\in\N$, $n>1$ and  $f\in\Pn$, we wish  to study the family
\begin{equation}
D_n(f)=\{g\in \PD: \quad  g^n=f\}
\end{equation}
 and  the quantity $\car( D_n(f))$, i.e.,  the cardinality of $D_n(f)$.
Note that    there exists $f\in\Pn\setminus\Pe$ such that the factor $g$ in (3)  will generally not unique. It turns out that $\car( D_n(f))$ depends on $n$ and in  some way also depends on the  geometric structure of the zeros set $N_f=\{x\in\R: f(x)=0\}$  of  $f\in\Pn$ (see our Theorem 1 below). The essential support $S_f$ of $f\in\PD$ is defined by $S_f=\R\setminus N_f$.
Combining (1) with (2), gives
\begin{equation}
 0\in S_f  \qquad {\text{and}}\qquad -S_f=S_f .
\end{equation}
Since functions $f\in\PD$ are continuous on $\R$, it follows that $S_f$ is an open subset of $\R$. Therefore,  $S_f$   can  be represented as a finite or infinite  union $S_f=\bigcup_{j\in \Sigma}E_j$, where $\{E_j\}_{j\in \Sigma}$ is the family of all  open connected components  of  $S_f$.  In the sequel,    ${\text{comp}}\,(S_f)$  denotes the cardinality of $\Sigma$. According to (5), we see that either there is an $k\in\N$ such that ${\text{comp}}\,(S_f)=2k-1$ or ${\text{comp}}\,(S_f)=\infty$.
\begin{thm}
 \ Let $n\in\N$, $n\ge 2$,  and let $f\in\Pn$. Assume that
\begin{equation}
{\text{comp}}\,(S_f)=2k-1
\end{equation}
for some  $k\in\N$. Then
\begin{equation}
\car(D_n(f))\le n^{k-1}.
\end{equation}
\end{thm}
The following theorem shows that the estimate (7) is accurate.
\begin{thm}
 \ For each $ n\in\N$, $n\ge 2$, each $k\in\N$,  and any  open subset $E$ of $ \R$ which satisfies
\begin{equation}
0\in E, \qquad -E=E \qquad{\text{and}}\qquad {\text{comp}}\,(E)=2k-1,
\end{equation}
 there exists $f\in \Pn$ such that $S_f=E$ and
\begin{equation}
 \car(D_n(f))= n^{k-1}.
\end{equation}
\end{thm}
We will present two examples of $f\in\Pn$ such that  (9) is satisfied. In order to make the examples easier to understand, we will  consider only small values of $k$. For $\alpha>0$, set
\begin{equation}
\Lambda_{\alpha}(x)=\max\{(1-|x|)^{\alpha}; \ 0\},
\end{equation}
 $x\in\R$.  Note that    $\Lambda_{\alpha}\in \PD$ if and only if $\alpha\ge 1$ (see, e.g., \cite[p. 282]{16}). We start with the case where $S_f$ is a bounded subset  of $\R$.
\begin{exam}
For $ n\in\N$, $n>1$,  let
\begin{gather}
f_1(x)=\Lambda_n(x)\nonumber\\
+\frac1{2^{3n}}\biggl[\Bigl(\Lambda_1(x-\pi-1)+\Lambda_1(x-2\pi+1)\Bigr)^n+\Bigl(\Lambda_1(x+\pi+1)+ \Lambda_1(x+2\pi -1)\Bigr)^n\biggr]\nonumber\\
+\frac1{2^{4n}}\biggl[\Bigl((\Lambda_1(x-11)+\Lambda_1(x-12)+\Lambda_1(x-13)+\Lambda_1(x-14)\Bigr)^n\nonumber\\
+\Bigl(\Lambda_1(x+11)+\Lambda_1(x+12)+ \Lambda_1(x+13)+\Lambda_1(x+14)\Bigr)^n\biggr].
\end{gather}
Then ${\text{comp}}\,(S_{f_1})=5$, since
\begin{equation}
S_{f_1}=(-15,-10)\cup (-2\pi,-\pi)\cup(-1,1)\cup (\pi, 2\pi)\cup(10,15).
\end{equation}
Moreover,  $f_1\in PD_n(\R)$ and
\begin{equation}
 \car(D_n(f_1))= n^5.
\end{equation}
\end{exam}
Now let us give an example of $ f\in \Pn$ such that $S_{f}$ has an unbounded components.
\begin{exam}
For any $a>15$ and $ n\in\N$, $n>1$,  let
\begin{gather}
f_2(x)=f_1(x)+\frac1{2^{2n}}\biggl(\sum_{r=1}^{\infty}\frac1{2^{r}} \Lambda_1(x-a-r)\biggr)^n+\biggl(\sum_{r=1}^{\infty}\frac1{2^{r}} \Lambda_1(x+a+r)\biggr)^n,
\end{gather}
where $f_1$ was defined by (11). Then ${\text{comp}}\,(S_{f_2})=7$, since
\begin{equation}
S_{f_2}=(-\infty, -a)\cup(-15,-10)\cup (-2\pi,-\pi)\cup(-1,1)\cup (\pi, 2\pi)\cup(10,15)\cup(a,\infty).
\end{equation}
Moreover,  $f_2\in PD_n(\R)$ and
\begin{equation}
 \car(D_n(f_2))= n^7.
 \end{equation}
\end{exam}
\begin{thm}
Suppose that an open subset $E$ of $ \R$ satisfies
\begin{equation}
0\in E, \qquad -E=E \qquad{\text{and}}\qquad {\text{comp}}\,(E)=\infty.
\end{equation}
Let $\{E_j\}_{j\in \Sigma} $ be the family of open connected components  of $E$. Assume that there is $\sigma>0$ such that
\begin{equation}
\inf_{j\in \Sigma}\bigl(\sup_{a, b\in E_j} |a-b|\big)=2\sigma.
\end{equation}
Then,  for any $n\in\N$, $n\ge 2$,   there exists $f\in \Pn$ such that $S_f=E$ and
\begin{equation}
 \car(D_n(f))= \infty.
\end{equation}
\end{thm}

\section{Preliminaries and  Proofs}

 If $\nu, \mu\in \M$, then the convolution $\nu\ast\mu$ is defined by
 \[
 \nu\ast\mu(E)=\int_{\R}\nu(E-x)\,d\mu(x)
 \]
 for each Borel subset $E$ of $\R$. Note that
\begin{equation}
\widehat{\nu\ast\mu}=\widehat{\nu}\cdot\widehat{\mu}.
\end{equation}
In particular, for any $n\in\N$,
\begin{equation}
\widehat{\Bigl(\mu^{\ast n}\Bigr)}=\Bigl(\widehat{\mu}\Bigr)^{n},
\end{equation}
where the convolution power $\mu^{\ast n}$ is defined as the $n$-fold iteration of the convolution of $\mu$ with itself.

 The Lebesgue space $\Le$  can be identified with the closed ideal in $\M$ of   measures absolutely continuous with respect to the Lebesgue measure $dx$ on $\R$. Namely, if $\varphi\in \Le$, then  $\varphi$ is associated with  the  measure
\[
\mu_{\varphi}(E)=\int_{E}\varphi(t)\,dt
\]
for  each Borel subset $E$ of $\R$. Hence  $\widehat{\varphi}(x)=\int_{\R}e^{-itx}\varphi(t)\,dt$.  In particular, if $\mu=\varphi(t)dt$, where  $\varphi\in L^1(\R)$ and $\varphi$ is such  that  $\|\varphi\|_{L^1(\R)}=1$ and $\varphi\ge 0$ on $\R$,   then $\varphi$  is called the probability density function of $\mu$, or the probability density for short.

We  define the inverse Fourier transform by
 \[
\widecheck{\psi}(t)=\frac1{2\pi}\int_{\R}e^{it\,x}\psi(x)\,dx,
\]
$t\in\R$. Then  the inversion  formula  $\widehat{(\check{\psi})}=\psi$ holds for suitable $\psi\in L^1(\R)$.

{\bf{Proof of Theorem 1.}}\quad The conditions (5) and  (6) imply that there exists a sequence of real numbers
\begin{equation}
0<b_0<a_1<b_1<a_2<b_2<\dots<a_{k-1}<b_{k-1}
\end{equation}
such that
\begin{equation}
S_f=\Bigl(\bigcup_{j=-1}^{-(k-1)}E_j\Bigr)\bigcup E_0\bigcup\Bigl(\bigcup_{j=1}^{k-1}E_j\Bigr),
\end{equation}
where $E_0=(-b_0, b_0)$, $E_j=(a_j,b_j)$ and $E_{-j}=-E_j=(-b_j,-a_j)$ for $j=1,2,\dots, k-1$. Note that  in (22)   also might be $b_{k-1}=+\infty$.
Let $g\in\Dn$. Then it is immediate that $S_f=S_g$ and
\begin{equation}
|g(x)|^n=|f(x)|
\end{equation}
for all $x\in S_f$. Fix any $E_j\subset S_f$ in (23). Since $E_j$ is an open connected component of $S_f$, we have that there are two continuous functions $u_{f,j}, u_{g,j}: E_j\to (-\pi, \pi]$ such that
\begin{equation}
f(x)=|f(x)|e^{iu_{f,j}(x)}\quad {\text{and}}\quad  g(x)=|g(x)|e^{iu_{g,j}(x)}
\end{equation}
for all   $x\in E_j$. Using the identity  $ g^n=f$, it follows from (24) and (25) that, for each  for  $j\in \{-(k-1),\dots, k-1\}$, there exists some integer $m_j$ in $\{0,\dots, n-1\}$ such that
\begin{equation}
n\cdot u_{g,j}(x)= u_{f,j}(x)+2\pi m_j
\end{equation}
for all $x\in E_j$.  Therefore, for any $x\in S_f$, we have
\begin{equation}
g(x)=\sum _{j=-(k-1)}^{k-1}|g(x)|\chi_{E_{j}}(x) e^{iu_{g,j}}(x)=\sum _{j=-(k-1)}^{k-1}|f(x)|_{+}^{1/n}\chi_{E_{j}}(x) e^{(iu_{f,j}+2\pi m_j)/n}(x),
\end{equation}
where $\chi_{E_{j}}$ is the indicator function of the set $E_j$  and $|f(x)|_{+}^{1/n}$  denotes the positive $n$th root of positive number $|f(x)|$, $x\in S_f$.  We claim that $m_0=0$. Indeed, (2) implies that $f(0)$ and $g(0)$ are positive numbers.  Then (25) implies that $u_{f,0}(0)=u_{g,0}(0)=0$.
Combining this with (26), yields the claim. Next, applying the property (1), we get that $m_{-j}=-m_j$  for all $j\in\{0,1,\dots, k-1\}$.
Therefore, we conclude from  (27) that
\begin{gather}
g(x)=|f(x)|^{1/n}_{+}\chi_{E_{0}}(x) +\sum _{j=1}^{k-1}|f(x)|^{1/n}_{+}\Bigl( \chi_{E_{j}}(-x) e^{-(iu_{f,j}(x)+2\pi m_j)/n}\nonumber \\
+\chi_{E_{j}}(x) e^{(iu_{f,j}(x)+2\pi m_j)/n}\Bigr)
\end{gather}
for all $x\in S_g$. Finally, keeping in mind that each $m_j$, $j=1,2,\dots,k-1$, may take any  value in $\{0,1,\dots,n-1\}$, we obtain from (28) the estimate (7).   Theorem 1 is proved.
\begin{rem}
 {\rm{Of course, we are not claiming   that each of  $n^{k-1}$th possible functions in (28)  belongs to $\PD$.}}
 \end{rem}
\begin{rem}
 {\rm{In the proof of Theorem 1 we concerned with the so-called problem of phase retrieval (see the  equality  (24)), i.e., the problem of the recovery of a measure $\mu$  given the amplitude $|f|$ of its Fourier transform $f=\widehat{\mu}$. This problem is well known in various fields of science and engineering, including crystallography, nuclear magnetic resonance and
optics (see, for example, survey [4]).}}
 \end{rem}
 {\bf{Proof of Theorem 2.}}\quad
 Note that as in the proof of Theorem 1, in light of (8),  we see  that  there exists a sequence of real numbers (22) such that
\begin{equation}
E=\Bigl(\bigcup_{j=-1}^{-(k-1)}E_j\Bigr)\bigcup E_0\bigcup\Bigl(\bigcup_{j=1}^{k-1}E_j\Bigr),
\end{equation}
Let us split our proof into two cases.
First  we consider the case when all $E_j$ in (29) are finite intervals. Denote by $2\sigma$ the minimal length of $E_j$, $j=-(k-1),\dots, k-1$, i.e.,
\begin{equation}
2\sigma= \min\Bigl\{2b_0; \min_{1\le j\le k-1}(b_j-a_j)\Bigr\}.
\end{equation}
Let $\varphi\in \PD$. Assume, in addition,  that $\varphi$  is real-valued on $\R$  and
\begin{equation}
S_{\varphi}=(-\sigma, \sigma).
\end{equation}
For example, we can take $\varphi(x)=\Lambda_{\alpha}(x/\sigma)$, where the truncated power function $\Lambda_{\alpha}$ was defined  by (9).
 Next, for each $j\in\{1,2,\dots , k-1\}$,  we take any sequence of real numbers $\{\tau_{j,s}\}_{s=1}^{m(j)}$ such that
\begin{equation}
a_j+\sigma=\tau_{j, 1}<\tau_{j, 2}<\dots <\tau_{j, m}=b_j-\sigma.
\end{equation}
In addition, we assume that
\begin{equation}
 \tau_{j, (s+1)}- \tau_{j, s}<2\sigma
\end{equation}
for all $s=1,\dots, m(j)$. Then we define the function
\begin{equation}
u_j(x)=\sum_{s=1}^{m(j)}\varphi(x-\tau_{j, s}),
\end{equation}
$x\in\R$. Now  (32) and (33) imply that $u_j$ is supported on $[a_j,b_j]=\overline{E_j}$ and $S_{u_{j}}= E_j$.  Since $\varphi\in \PD$ is real-valued on $\R$ and satisfies (31),  we conclude that $\varphi$  is  even and positive  on $S_{\varphi}$. Therefore, the function (34) and the function
\begin{equation}
u_j(-x)=\sum_{s=1}^{m(j)}\varphi(x+\tau_{j, s})
\end{equation}
are strictly positive on $E_j$ and on $E_{-j}=-E_j$, respectively.

Let us   define the function $u_{0}$ such that it  is supported on $[-b_0, b_0]=\overline{E{_0}}$ and $S_{u_{0}}= E_{0}$. To this end, we take any sequence $\{\theta_{i}\}_{i=
0}^{r}$ of real numbers such that
\begin{equation}
-b_{0}+\sigma=\theta_{r}<\dots \theta_{1}<\theta_0=0
\end{equation}
and
\begin{equation}
 \theta_{i}- \theta_{i+1}<2\sigma
\end{equation}
for $i=0,\dots, l-1$. Next,  for an arbitrary  sequence of positive numbers $\{\omega_i\}_{i=0}^{l}$, we define
\begin{equation}
u_0(x)=\sum_{i=0}^{l} \omega_i\Bigl(\varphi(x-\theta_i)+\varphi(x+\theta_i)\Bigr),
\end{equation}
$x\in \R$. Of course, (36) and (37) imply that $u_{0}$  is supported on $[-b_0, b_0]=\overline{E{_0}}$ and $S_{u_{0}}= E_{0}$.

Finally, given any fixed    sequence  of  positive numbers numbers  $\{\alpha_j\}_{j=1}^{k-1}$,  we set
\begin{equation}
u(x)=u_{0}(x)+\sum_{j=1}^{k-1}\alpha_j\Bigl(u_j(x)+u_j(-x)\Bigr),
\end{equation}
$x\in \R$.

We claim that $u\in\PD$ and $S_u=E$.
First, as real-valued function $\varphi\in\PD$ satisfies (31), it follows from  (29)-(33) and from, (36)-(37) that $u$ is  continuous  on $\R$ and
\[
S_u=\bigcup_{j=-k}^{k}E_j =E.
\]
Second, since $u$ is continuous and compactly supported, it follows that $u\in L^1(\R)$. Therefore, the inverse Fourier transform of $u$ is well-defined. Hence
\begin{gather}
\widecheck{u}(t)=\widecheck{u_0}(t)+\sum_{j=1}^{k-1}\Bigl(\alpha_j\widecheck{\varphi}(t)\sum_{s=1}^{m(j)}2\cos(\tau_{j,s}\cdot t)\Bigr)\nonumber \\
=2\widecheck{\varphi}(t)\Biggl[\sum_{i=0}^l\omega_i\cos(\theta_i\cdot t)+\sum_{j=1}^{k-1}\Bigl(\alpha_j\sum_{s=1}^{m(j)}\cos(\tau_{j,s}\cdot t)\Bigr)\Biggr]\nonumber \\
=2\widecheck{\varphi}(t)\Biggl[\omega_0+\sum_{i=1}^l\omega_i\cos(\theta_i\cdot t)+\sum_{j=1}^{k-1}\Bigl(\alpha_j\sum_{s=1}^{m(j)}\cos(\tau_{j,s}\cdot t)\Bigr)\Biggr].
\end{gather}
 Let us fix the previously chosen positive numbers $\omega_1,\dots, \omega_k$ and $\alpha_1,\dots,\alpha_{k-1}$.  Then we  increase if necessary,  the value of $\omega_0$  in such a way that
\begin{equation}
\omega_0>\sum_{i=1}^l\omega_i  + \sum_{j=1}^{k-1}\alpha_j m_j.
\end{equation}
 Bochner's theorem shows that $\widecheck{\varphi}(t)\ge 0$ for all $t\in\R$, since $\varphi\in \PD$. Combining (31) with (41),  we see that $\widecheck{u}$ is nonnegative on $\R$ and $\widecheck{u}\not\equiv 0$.  In addition, we conclude (see, e.g., \cite[p. 409]{1}) that  $\widecheck{u}\in L^1(\R)$. Thus, applying  the  Fourier transform to $\widecheck{u}$ and  using again Bohner's theorem, we see  that $u$ is continuous nontrivial positive definite, i.e., $u\in \PD$. This proves our claim.

Define
\begin{equation}
f(x)=u_{0}^n(x)+\sum_{j=1}^{k-1}\alpha_j^n\Bigl(u_j^n(x)+u_j^n(-x)\Bigr),
\end{equation}
$x\in\R$. We claim that $f$ satisfies the hypotheses of Theorem 2. Let us first prove that $f\in\PD$. Indeed, from (22) we see that the essential support $S_u$ of $u$, defined by (39),  can be represented as the union $\bigcup_{j=-(k-1)}^{k-1} E_j$ of a family of pairwise disjoint sets $E_j=S_{u_j}$, $j=-(k-1)\dots,k-1$. Therefore,
\begin{equation}
u^n(x)=u_{0}^n(x)+\sum_{j=1}^{k-1}\alpha_j^n\Bigl(u_j^n(x)+u_j^n(-x)\Bigr),
\end{equation}
$x\in\R$. We have already proven that $u\in\PD$. On the other hand, it is well known that for each $n\in\N$ and any $\zeta\in\PD$, it follows that $\chi^n\in\PD$. Thus, $u^n\in\PD$. Combining this fact with (36), we conclude that $f\in\PD$ and $S_f=S_u=\bigcup_{j=-(k-1)}^{k-1}E_j=E$.

Second,  we will prove that the function $f$ defined by (43) has the property (9). To this end, let $\Zn$ denote the group $\Z_n=\Z\diagup n\Z\cong \{0,1,\dots, n-1\}$. Given  $\Lambda=(p_1,p_2,\dots, p_{k-1})\in \Zn^{k-1}$,   define
\begin{equation}
g_{\Lambda}(x)=u_0(x)+\sum_{j=1}^{k-1}\alpha_j\biggl(e^{i2\pi p_j/n}u_j(x)+e^{-i2\pi p_j/n}u_j(-x)\biggr),
\end{equation}
$x\in\R$.
 We claim that $g_{\Lambda}\in D_n(f)$. By the same argument as before for the function $u$  defined by (39), we see that  $g_{\Lambda}\in L^1(\R)$. Therefore, $\widecheck{g_{\Lambda}}$ is well-defined and
\begin{equation}
\widecheck{g_{\Lambda}}(t)=2\widecheck{\varphi}(t)\biggl[\omega_0+\sum_{i=1}^l\omega_i\cos(\theta_i\cdot t)+\sum_{j=1}^{k-1}\biggl(\alpha_j\sum_{s=1}^{m(j)}\cos\biggl(\tau_{j,s}\cdot t+\frac{2\pi p_j}{n}\biggr)\biggr)\biggr],
\end{equation}
$t\in\R$.  Now using (41), we get that $\widecheck{g_{\Lambda}}(t)\ge 0$ for all $t\in\R$ and $ \widecheck{g_{\Lambda}}\not\equiv 0$.  Hence, by Bochner's theorem it follows that  $g_{\Lambda}\in \PD$. Next,
\begin{gather}
g_{\Lambda}^n(x)=u_{0}^n(x)+\sum_{j=1}^{k-1}\alpha_j^n\biggl(\biggl(e^{i2\pi p_j/n}u_j(x)\biggr)^n+\biggl(e^{-i2\pi p_j/n}u_j(-x)\biggr)^n\biggr)\nonumber \\
=u_{0}^n(x)+\sum_{j=1}^{k-1}\alpha_j^n\biggl(u_j(x)^n+u_j(-x)^n\biggr),
\end{gather}
$x\in\R$. Combining this representation with (42) and (43), we see that $g_{\Lambda}\in D_n(f)$, which yields our claim.

Finally, again using the fact that   $S_{g_{\Lambda}}$  is the union of a family of pairwise disjoint sets $S_{u_j}$, $j=-(k-1),\dots, k-1$, we  conclude from (44) that $g_{\Lambda}\equiv g_{\Lambda_1}$ for some $\Lambda, \Lambda_1\in \Zn^{k-1}$,  if and only if  $\Lambda=\Lambda_1$. This proves (9) in the case where each $E_j$ in (28) is a finite interval.

Now we consider the second case with $b_{k-1}=\infty$, i.e., if in (28) we have $E_{k-1}=(a_{k-1},\infty)$. Using the same $\varphi\in PD$ satisfying (31), we define the functions $u_0$ and $u_j$, $j=-(k-2),\dots, -1, 1,\dots,  k-2$ by (38) and (34)-(35), respectively. For $j=k-1$, let us take an arbitrary sequence of  positive numbers $\{\gamma_r\}_{r=1}^{\infty}$ such that
\begin{equation}
\sum_{r=1}^{\infty} \gamma_r=1.
\end{equation}
Then we define
\begin{equation}
u_{k-1}(x)=\sum_{r=1}^{\infty} \gamma_r\varphi(x-a_k-r\sigma),
\end{equation}
$x\in\R$. Obviously, ${\text{supp}}( u_{k-1})=\overline{E_{k-1}}$ and    $S_{u_{k-1}}=E_{k-1}$. Next, for the function $u$, defined by (39), it follows from (40) that
\begin{gather}
\widecheck{u}(t)=2\widecheck{\varphi}(t)\biggl[\omega_0+\sum_{i=1}^l\omega_i\cos(\theta_i\cdot t)+\sum_{j=1}^{k-2}\Bigl(\alpha_j\sum_{s=1}^{m(j)}\cos(\tau_{j,s}\cdot t)\Bigr)\nonumber \\
+\alpha_{k-1}\sum_{r=1}^{\infty}\gamma_r\cos\Bigl((a_{k-1}+r\sigma)t\Big)\biggr].
\end{gather}
Again, for fixed positive numbers $\omega_1,\dots \omega_k$, $\alpha_1,\dots, \alpha_{k-1}$, and $\gamma_1, \gamma_2,\dots$, we take the value of $\omega_0$ in such a way that
\begin{equation}
\omega_0>\sum_{i=1}^l\omega_i  + \sum_{j=1}^{k-2}\alpha_j m_j +\alpha_{k-1}.
\end{equation}
Combining (31) with (50), we conclude  from (49) that $\widecheck{u}$ is nonnegative on $\R$ and $\widecheck{u}\not\equiv 0$. Thus,  $u\in\PD$. Finally,  we claim that the function $f$ defined by (42) also satisfies the hypotheses of Theorem 2 in our case with $E_{k-1}=(a_{k-1}, \infty)$. The  proof of this claim     is exactly the same as that of  the first case. Therefore, we skip the details of this proof. Theorem 2 is proved.

{\bf{Proof of Example 3.}}\quad We claim that there are $u_0$, $u_1$, $u_2$ and $\alpha_1$, $\alpha_2$ such that the function in (11) coincides with the function  defined by (43). Indeed, set $\varphi=\Lambda_1$, $j=3$, $m_1=2$, $m_2=4$ and
\[
\tau_{11}=\pi+1, \ \tau_{12}=2\pi-1, \ \tau_{21}=11, \ \tau_{22}=12, \ \tau_{23}=13, \ \tau_{24}=14.
\]
Since $\sigma$ defined by (30) is equal now to $1$, then $\tau_{js}$ defined above satisfies (32) and (33). Next, set
\[
 v_1(x)=\Lambda_1(x-\pi-1)+\Lambda_1(x-2\pi+1)
 \]
 and
 \[
 v_2(x)=\Lambda_1(x-11)+\Lambda_1(x-12)+\Lambda_1(x-13)+\Lambda_1(x-14).
 \]
 For $v_0=\Lambda_1$,  $\alpha_1=1/8$ and $\alpha_2=1/16$, let $v$ be given by
 \[
 v(x)=\Lambda_{1}(x)+ \frac18\Bigl(v_1(x) +v_1(-x)\Bigr)+\frac1{16}\Bigl(v_2(x)+v_2(-x)\Bigr).
 \]
 It is easily seen that $\Lambda_1$, $v_1$ and $v_2$ are supported on a family of pairwise disjoint sets. Therefore,
 \begin{equation}
 v^n(x)=\Lambda_{n}(x)+ \frac18\Bigl(v_1^n(x) +v_1^n(-x)\Bigr)+\frac1{16}\Bigl(v^n_2(x)+v^n_2(-x)\Bigr),
  \end{equation}
 since $\Lambda^n_1=\Lambda_n$. Next,
 \begin{equation}
 S_{v^n}= S_v=(-15,-10)\cup (-2\pi,-\pi)\cup(-1,1)\cup (\pi, 2\pi)\cup(10,15).
 \end{equation}
   The function $v^n$ coincides with the function in (11) and is defined by the same rules as in (43).  Our claim is proved.

  Now, it is enough to show that $v\in\PD$. Indeed,
 \begin{gather}
\widecheck{v}(t)=2\widecheck{\Lambda}_1(t)\biggl[1+\sum_{j=1}^{2}\Bigl(\alpha_j\sum_{s=1}^{m(j)}\cos(\tau_{js}\cdot t)\Bigr)\biggr]\nonumber \\
=2\widecheck{\Lambda}_1(t)\biggl[1+\frac18\Bigl(\cos(\tau_{11}t)+\cos(\tau_{12}t)\Bigr)\nonumber \\
+\frac1{16}\Bigl(\cos(\tau_{21}t)+\cos(\tau_{22}t)+\cos(\tau_{23}t)+\cos(\tau_{24}t)
\Bigr)\biggr]\ge\widecheck{\Lambda}_1(t)\ge 0
\end{gather}
 for all $t\in\R$, since $\Lambda_1\in\PD$. Therefore, Bochner's theorem shows that $v\in\PD$.

 By repeating the finally  part of the proof of Teorema 2, we complete the proof  of Example 3.

{\bf{Proof of Example 4.}}\quad This example concerns the case that was considered in the second part of the proof of  Theorem 2, i.e., when $S_{f}$ contains two unbounded components $(-\infty, -a)$ and $(a, \infty)$. Also, as in the proof of Theorem 2, is enough to show that (47) and (50) are satisfied. Indeed, (47) is clear, since $\gamma_r=1/2^r$, $r=1,2,\dots$. We conclude from (11) and (14) that 
\[
k=4, \ \sigma=1, \ m_1=2,\ m_2=4,\ \omega_0=1, \omega_1=\cdot =0, \alpha_1=\frac18, \  \alpha_1=\frac1{16} \ {\text{ and}} \  \alpha_3=\frac1{4}.
\]
Therefore, a simple calculation shows that (50) is also satisfied. This completes the proof.

{\bf{Proof of Theorem 5.}}\quad We will prove this theorem using essential the same techniques as in the proof of Theorem 2. Therefore,   we sketch the proof only.  From (17) and (18) we have that there exits  an infinite sequence
\begin{equation}
0<b_0<a_1<b_1<a_2<b_2<\dots <\infty
\end{equation}
such that
\begin{equation}
E=\bigcup_{j\in\Z} E_j,
\end{equation}
where  $E_0=(-b_0, b_0)$ and $ E_j=(a_j, b_j)=-E_{-j}$ for each $j\in \N$. Moreover, from (18) we see that
\begin{equation}
\min\Bigl\{2b_0;\quad \inf_{j\in\N}(b_j-a_j)\Bigr\}=2\sigma>0.
\end{equation}
Let $\varphi\in\PD$ be the same function satisfying (31). For $j=0$ and for $j\in\N$, we define the functions $u_0$ and $u_j$ by (34) and (38). respectively.  Note that the sequences $\{\tau_{j_s}\}_{s=1}^{m(j)}$ and $\{\theta_i\}_{i=0}^l$ satisfy (32)-(33) and (35)-(37), respectively. For any $j\in\N$, let us define
\begin{equation}
\alpha_j=  \frac1{2^{j}m(j)\bigl|E_j\bigr|},
\end{equation}
where $|E_j|$ is the length of $E_j$, i.e., $|E_j|=\sup\{|a-b|: \ a,b\in E_j\}$. Note that the condition (18) guarantees that $\{\alpha_j\}_{j\in\N}$ is a well-defined sequence of positive numbers. Set
\begin{equation}
u(x)=u_{0}(x)+\sum_{j\in \N}\alpha_j\Bigl(u_j(x)+u_j(-x)\Bigr),
\end{equation}
$x\in \R$.  Combinnig (2), (34), (38) with (57), we get
\begin{equation}
|u(x)|\le 2\varphi(0)\biggl[\Bigl(\sum_{i=0}^l\omega_i\Bigr)I_{E_{0}}(x)+\sum_{j\in\N}\Bigl(\frac1{2^j|E_j|}I_{E_{j}}(x)\Bigr)\biggl]
\end{equation}
 for any $x\in\R$. Here, $I_D(x)$ denotes the indicator function of a subset $D\subset \R$. Hence,
 \[
 \|\varphi\|_{L^1(\R)}\le 2\varphi(0)\biggl[2b_0\sum_{i=0}^l\omega_i+1\biggr].
 \]
Therefore,  $\varphi\in L^1(\R)$ and
\begin{gather}
\widecheck{u}(t)=2\widecheck{\varphi}(t)\Biggl[\omega_0+\sum_{i=1}^l\omega_i\cos(\theta_i\cdot t)+\sum_{j\in\N}\Bigl(\alpha_j\sum_{s=1}^{m(j)}\cos(\tau_{j_{s}}\cdot t)\Bigr)\Biggr].
\end{gather}
Again, for fixed positive numbers $\omega_1,\dots \omega_k$, we take the value of $\omega_0$ in such a way that
\begin{equation}
\omega_0>\sum_{i=1}^l\omega_i  + \frac1{2\sigma}.
\end{equation}
Combining this condition with (56 and (57), we conclude   that the function $\widecheck{u}$ in (51) is nonnegative on $\R$ and $\widecheck{u}\not\equiv 0$. Thus,  $u\in\PD$.
Define
\begin{equation}
f(x)=u_{0}^n(x)+\sum_{j\in\N}\alpha_j^n\Bigl(u_j^n(x)+u_j^n(-x)\Bigr),
\end{equation}
$x\in\R$. Next, for $\Lambda=(p_1, p_2,\dots )\in \Z^{\infty}$, set
\begin{equation}
g_{\Lambda}(x)=u_0(x)+\sum_{j\in\N}\alpha_j\biggl(e^{i2\pi p_j/n}u_j(x)+e^{-i2\pi p_j/n}u_j(-x)\biggr),
\end{equation}
$x\in\R$. We claim that: \  (i)\  $f\in \PD$; (ii)  $E$ defined by(55) satisfies (17) and (18); \  (iii) \ $g_{\Lambda}\in D_n(f)$. \
The proofs of these  claims is exactly the same as in the case of functions $f$ and $g_{\Lambda}$ defined by (42) and (44), respectively.

Finally, again using the fact that $S_f=S_u$ and therefore,  $S_{g_{\Lambda}}=S_u$  are the unit of a family of pairwise disjoint sets $S_{u_j}$, $j\in\Z$, where $u_{-j}(x)=u_j(-x)$, $j\in\N$, we see  from (63) that $g_{\Lambda}\equiv g_{\Lambda_1}$ for some $\Lambda, \Lambda_1\in \Z^{\infty}$,  if and only if  $\Lambda=\Lambda_1$. This proves (19) and therefore completes the proof our theorem.

\section{Conclusion}

We study the $n$-divisible functions in $\PD$, where $\PD$ denotes the family of continuous positive definite functions on the real line $\R$. While there is a rich literature on infinite-divisible functions in $\PD$, for an integer $n>1$,  properties of $n$-divisible functions from $\PD$ have been studied much less. Surprisingly, it appears that  some properties of infinite-divisible and $n$-divisible functions $f\in\PD$  may differ in essence.  In this paper, we examine one such property, which has not yet been discussed in detail in the literature. More precisely, if $f\in\PD$ infinite-divisible, then   it is well-known that,  for each integer $n>1$, there is an unique $g\in\PD$, such that $g^n=f$. On the other hand, there is  $n$-divisible $f$ such that the factor $g$ in $g^n=f$ is generally not unique. For $n$-divisible $f\in\PD$, we study the following questions: \ (i) how rich can be the class $D_n(f)=\{g\in\PD: g^n=f\}$; \ (ii) what properties of $f$  determine the size of  $D_n(f)$. We present several
precise estimates for the cardinality of $D_n(f)$. Also, the main results are validated via illustrative examples.

\vspace{5mm}

{\bf {\Large{ Conflicts of Interest}}}

\vspace{3mm}

The author declare that they have no competing interests.

}}
\end{document}